\documentclass{amsart}

\usepackage{amscd}
\usepackage{amssymb}

\newtheorem{theo}{Theorem}
\newtheorem{prop}[theo]{Proposition}

\newtheorem{rema}[theo]{Remark}
\newtheorem{ex}[theo]{Example}

\title{Cosemisimple Hopf Algebras with Antipode of Arbitrary Finite Order}
\date{} 
\author{Julien Bichon}

\begin{document}

\makeatletter
\renewcommand{\@makefnmark}{}
\makeatother

\begin{abstract}
Let $m \geq 1$ be a positive integer. We show that, over
an algebraically closed field of characteristic zero,
there exists cosemisimple Hopf algebras having an antipode of order
$2m$. We also discuss the Schur indicator for such Hopf algebras.
\end{abstract}

\maketitle


\section{Introduction}

The order of the antipode of a Hopf algebra has been an important subject
of study since the beginning of the theory in the 60's.
Let us first recall some of the highlights.

\noindent
$\bullet$ First, there is the following fact: the antipode of 
a commutative or cocommutative Hopf algebra is involutive.

\noindent
$\bullet$ In 1971, Taft \cite{[Ta]} has constructed finite-dimensional
Hopf algebras with antipode of arbitrary even order. The Taft
algebras are not cosemisimple.

\noindent
$\bullet$ In 1975, Kaplansky \cite{[K]} conjectured that the antipode of 
 finite-dimensional cosemi-simple Hopf algebra is involutive.

\noindent
$\bullet$ Radford (\cite{[R]}, 1976) proved that the order
of the antipode of a finite-dimensional Hopf algebra is finite.

\noindent
$\bullet$
Kaplansky's  conjecture was proved in characteristic zero by Larson-Radford 
(\cite{[LR1],[LR2]}, 1987).
The conjecture remains opened in general,
although Etingof-Gelaki (\cite{[EG1]},
1998) proved it in positive characteristic,
under the additional assumption of semisimplicity.

\medskip

It seems that the only known examples of  cosemisimple
Hopf algebras have either involutive antipode or antipode
of infinite order
(the antipode of a cosemisimple Hopf algebra is always bijective). 
In this note, we present, under the assumption that
the base field is algebraically closed of characteristic zero, examples 
of cosemisimple Hopf algebras with antipode of arbitrary even order:

\begin{theo}
Let $m \geq 1$ be a positive integer and let $k$ be
an algebraically closed field of characteristic zero.
There exists a cosemisimple Hopf algebra over $k$
having an antipode of order
$2m$.
\end{theo}

In fact the Hopf algebras of the Theorem were introduced by Dubois-Violette
an Launer \cite{[DVL]}. Their properties are straightforward consequences 
of the results of \cite{[Bi]}. Using these Hopf algebras,
we also have a partial negative answer to a question raised in \cite{[EG2]}
(see Remark 5).

In the last Section, we study the Schur indicator for cosemisimple
Hopf algebras.
We get a generalization of the Frobenius-Schur theorem
of Linchenko-Montgomery \cite{[LM]} to the case of cosemisimple Hopf
algebras with involutive antipode. Finally we use the 
Hopf algebras of the second section to illustrate the difficulty
for proving a wider generalization, even when
the order of the antipode is finite.

\smallskip
 
We work over a fixed algebraically closed field of characteristic zero $k$.

\section{Proof of Theorem 1}

We consider the universal Hopf algebra associated to a non-degenerate
bilinear form, introduced by Dubois-Violette and Launer
\cite{[DVL]}. 
Let $n \in \mathbb N^*$,
and let $E \in GL_n(k)$. We consider the following algebra
$\mathcal B(E)$: it is the universal algebra with generators
$(a_{ij})_{1 \leq i,j \leq n}$ and satisfying the relations
$$E^{-\!1} {^t \! a} E a = I = a E^{-\!1} {^t \! a} E,$$
where $a$ is the matrix $(a_{ij})_{1 \leq i,j \leq n}$ and $I$ is the
identity matrix. The algebra $\mathcal B(E)$ admits a 
Hopf algebra structure, with
comultiplication $\Delta$ defined by $\Delta(a_{ij})
= \sum_{k=1}^n a_{ik} \otimes a_{kj}$, $1 \leq i,j \leq n$,
with counit $\varepsilon$ defined by $\varepsilon(a_{ij}) = \delta_{ij}$,
$1 \leq i,j \leq n$, and with antipode $S$ defined on the matrix 
$a = (a_{ij})$ by
$S(a) = E^{-1}{^t \! a} E$.

\smallskip

For $\alpha_1 , \ldots , \alpha_n \in k$, the corresponding anti-diagonal
matrix is denoted by \break
AD$(\alpha_1 , \ldots , \alpha_n)$
(the coefficient $\alpha_n$ is located at the top-right of the matrix).
The corresponding diagonal matrix is denoted by
D$(\alpha_1 , \ldots , \alpha_n)$.

\smallskip

Now fix $m \geq 1$ and $\xi \in k^*$ a primitive $m$-th root of unity.
We consider the matrices
$$
E = {\rm AD}(\xi, 1, 1 ,1, 1 , 1)
\quad {\rm and} \quad 
F = E^{-1} {^t \!E} = {\rm D}(\xi^{-1}, 1,1,1,1, \xi).$$

\begin{prop}
The Hopf algebra $\mathcal B(E)$ is cosemisimple,
with antipode of order $2m$.
\end{prop} 

\noindent
\textbf{Proof}.
We first prove that $\mathcal B(E)$ is cosemisimple. Let
$q \in k^*$ be such that
$q^2 + {\rm tr}(F)q +1 = 0$. Then by \cite{[Bi]}, Theorem 1.1, there exists 
an equivalence of monoidal categories
$${\rm Comod}(\mathcal B(E)) \cong^{\otimes} 
{\rm Comod}(\mathcal O(SL_q(2))).$$
It is well-known that $\mathcal O(SL_q(2))$ is cosemisimple if
and only if $q = \pm 1$ or if $q$ is not a root of unity
(see e.g. \cite{[KS]}).
For the value of $q$ chosen here, it is easily seen that
either $q = -1$ or $q$ is not a root of unity
(e.g. embedding $\mathbb Q(\xi)$ into $\mathbb C$). 
Hence $\mathcal B(E)$ is cosemisimple.

Since $S^2(a) = FaF^{-1}$, we have
$S^{2m}(a) =a$ and the antipode of $\mathcal B(E)$ has order $\leq 2m$.
Now consider the 6-dimensional comodule
associated to the elements $a_{ij}$: this comodule corresponds,
via the category equivalence of \cite{[Bi]},
to the simple 2-dimensional $\mathcal O(SL_q(2))$-comodule.
It follows that the elements $a_{ij}$, $1 \leq i,j \leq 6$, 
are linearly independent. 
Let $G \in M_6(k)$ be such that $aG = G{^t a}$. The linear independence
of the $a_{ij}$'s forces $G = 0$, and thus it is clear that
the antipode of $\mathcal B(E)$ has even order.
Now let $k \in \mathbb N^*$ be such that 
$S^{2k}(a) = F^k a F^{-k} =a$. Again by the linear independence of
the $a_{ij}$'s, there exists $\lambda \in k^*$ such that
$F^k = \lambda I$, and hence $\xi^k = 1$. Since $\xi$ is primitive
$m$-th root of unity, $m$ divides $k$ and we conclude that
the antipode of $\mathcal B(E)$ has order $2m$. $\square$

\medskip

Theorem 1 is an immediate consequence of Proposition 2.

\begin{rema}
{\rm The cosemisimplicity of $\mathcal B(E)$
cannot be proved using a compactness-like argument
(we assume in this remark that $k = \mathbb C$).
Indeed $\mathcal B(E)$ does not admit a CQG algebra structure
\cite{[KS]}. This is easily seen examining the eigenvalues of the matrix
$F$ (see \cite{[KS]}, Lemma 30 of Chapter 11).}
\end{rema}

\begin{rema}
{\rm The matrix $E$ just considered does not have the smallest possible
size for particular values of $m$. 
Here are the smallest sizes we have found. 
Again $\xi$ 
is  a primitive $m$-th root of unity.

\noindent
$\bullet$ Assume that $m \geq 5$. Put $E = {\rm AD}(1,1, \xi)$. Then 
$\mathcal B(E)$ is a cosemisimple Hopf algebra
with antipode of order $2m$.

\noindent
$\bullet$ Assume that $m=4$. Put $E = {\rm AD}(1,1 ,1, \xi)$.
Then $\mathcal B(E)$ is a cosemisimple Hopf algebra
with antipode of order $8$.

\noindent
$\bullet$ Assume that $m=3$. Put $E = {\rm AD}(1,1 ,\xi, \xi)$.
Then $\mathcal B(E)$ is a cosemisimple Hopf algebra
with antipode of order $6$, and is cotriangular since
the corresponding $q$ in Theorem 1.1 of \cite{[Bi]} is $q=1$.}
\end{rema}

\begin{rema}
{\rm Let $t \in k^*$ be such that $t^2 + 3t +1 = 0$. Let  
$E = {\rm AD}(1,1,t)$. The corresponding $q$ in Theorem 1.1 of \cite{[Bi]}
 is $q = 1$, and thus $\mathcal B(E)$ is a cotriangular Hopf
algebra. Since $t^2$ is not a root of unity and
 is an eigenvalue of 
$S^2$, we have a partial negative answer to Question 7.4 
in \cite{[EG2]} (of course $\mathcal B(E)$ is not the twist of a 
function algebra).}
\end{rema} 

\section{The Schur indicator for cosemisimple Hopf algebras}

Let $G$ be a compact group and let $V$ be a complex finite-dimensional
representation of $G$. Following the notation of
\cite{[LM]}, we define the Schur indicator of $V$ to be
$$\nu_2(V) = \int_G \chi_V(g^2)dg.$$
The classical Frobenius-Schur theorem (see \cite{[BD]})
 states that for an irreducible
representation $V$, then $\nu_2(V) = 0,1$ or $-1$, with 
$\nu_2(V) \not = 0$ if and only if $V$ is self-dual. The case
$\nu_2(V) = 1$ corresponds to the existence of a $G$-invariant
symmetric non-degenerate bilinear form on $V$, while
the case
$\nu_2(V) = -1$ corresponds to the existence of a $G$-invariant
skew-symmetric non-degenerate bilinear form on $V$.

\smallskip

The Frobenius-Schur theorem for finite groups was generalized to
finite-dimensional semisimple Hopf algebras by Linchenko-Montgomery
\cite{[LM]}. We prove such a theorem for cosemisimple Hopf algebras
with involutive antipode, and discuss the difficulty for
predicting a more general theorem, even for 
 cosemisimple Hopf algebras
with antipode of order 4.

\smallskip

Let $A$ be a cosemisimple Hopf algebra with Haar measure $h$, and let
$V$ be a finite-dimensional $A$-comodule with corresponding coalgebra
map $\Phi_V : V^* \otimes V \longrightarrow A$.
The character of $V$ \cite{[La]}
is defined to be $\chi_V := \Phi_V({\rm id}_V)$.
Dualizing \cite{[LM]}, we define the Schur indicator of $V$ to be 
$$\nu_2(V) := h(\chi_{V(1)} \chi_{V(2)}),$$ 
with Sweedler's notation $\Delta (a) = a_{(1)} \otimes a_{(2)}$.

\begin{theo}
Let $A$ be a cosemisimple Hopf algebra, and let
$V$ be a finite-dimensional irreducible $A$-comodule.

\noindent
1) If the $A$-comodule $V$ is not self-dual, then $\nu_2(V) = 0$.

Assume now that $V$ is self-dual.

\noindent
2) Let $\beta : V \otimes V \longrightarrow k$ be a non-degenerate
$A$-colinear bilinear form. Let $E$ be the matrix of $\beta$ in some basis
of $V$. Then 
$$\nu_2(V) = \frac{\dim (V)}{{\rm tr}(E {^t \! E}^{-1})}.$$ 
3) Assume that $S^2 = {\rm id}$. Then $\nu_2(V) = \pm 1$. 
 The case
$\nu_2(V) = 1$ corresponds to the existence of an $A$-colinear
symmetric non-degenerate bilinear form on $V$, while
the case
$\nu_2(V) = -1$ corresponds to the existence of an $A$-colinear
skew-symmetric non-degenerate bilinear form on $V$.
\end{theo}

\noindent
\textbf{Proof}. Let $e_1, \ldots ,e_n$ be a basis of $V$, and let
$a_{ij}$, $1 \leq i,j \leq n$, be the corresponding
matrix coefficients of the comodule $V$. Then 
$\chi_V = \sum_i a_{ii}$ and $\nu_2(V) = \sum_{i,j} h(a_{ji}a_{ij})$.
By the orthogonality relations \cite{[KS]}, Proposition 15 of Section 11
(the orthogonality relations first appeared in \cite{[La]}), we have
$\nu_2(V) = 0$ if $V$ is not self-dual. 

Assume now that $V$ is self-dual. Let $E$ be the matrix of 
$\beta : V \otimes V \longrightarrow k$ in the fixed basis of $V$.
Since $\beta$ is $A$-colinear, we have $S(a) = E^{-1} {^t \! a}E$ and
 $S^2(a) = E^{-1} {^t \! E} a {^t \! E}^{-1} E$. Then again by the 
orthogonality relations \cite{[KS]}, we have
$$h(a_{kl}S(a_{ij})) = \delta_{kj} 
\frac{(E^{-1} {^t \! E})_{il}}{{\rm tr}(E {^t \! E}^{-1})},
\ 1 \leq i,j,k,l \leq n.$$
Using $S(a) = E^{-1} {^t \! a} E$, a direct computation gives
$$h(a_{kl}a_{ij}) =  
\frac{E_{ki}^{-1} E_{lj}}{{\rm tr}(E {^t \! E}^{-1})},
\ 1 \leq i,j,k,l \leq n.$$
This leads to
$$\nu_2(V) = \sum_{i,j}
\frac{E_{ji}^{-1} E_{ij}}{{\rm tr}(E {^t \! E}^{-1})}
= \frac{\dim (V)}{{\rm tr}(E {^t \! E}^{-1})},$$
as claimed. Assume finally that $S^2 = {\rm id}$.
The linear independence of the $a_{ij}$'s forces 
$E^{-1}{^t \! E} = \lambda I$ for $\lambda  \in k^*$, and necessarily
$\lambda = \pm 1$. This gives $\nu_2(V) = \pm 1$.
The last assertion is immediate. $\square$

\medskip

Using the Larson-Radford's theorems \cite{[LR1],[LR2]}, 
it is clear that Theorem 6 implies Theorem 3.1 of \cite{[LM]}.

\smallskip

The Hopf algebras $\mathcal O(SL_q(2))$ show that the Schur indicator
may take various possible values.
Even in the case of cosemisimple Hopf algebras with antipode of
finite order $>2$, it seems difficult to predict all the possible
values of the Schur indicator. The following example might convince the 
reader.

\begin{ex}
{\rm Let $n \geq 3$, with $n$ odd.
We claim that there exists a cosemisimple Hopf algebra with antipode of
order $4$ having an irreducible comodule $V$ with
$\nu_2(V) =n$.

Put $k = \frac{n-1}{2}$. Consider the matrices
$$E = {\rm AD}(\underbrace{-1,-1,\ldots , -1}_{k \ times} , 
\underbrace{1,1, \ldots , 1}_{2n -k \ times}) \quad {\rm and}$$
$$E^{-1} {^t \! E} = 
{\rm D}(\underbrace{-1,-1,\ldots , -1}_{k \ times} , 
\underbrace{1,1, \ldots , 1}_{2n -2k \ times},
\underbrace{-1,-1,\ldots , -1}_{k \ times}).
$$
Then  ${\rm tr}(E {^t \! E}^{-1}) = 2n-4k = 2$, and hence
$\mathcal B(E)$ is cosemisimple by \cite{[Bi]}. Consider now the fundamental
$2n$-dimensional $\mathcal B(E)$-comodule
(irreducible by \cite{[Bi]}), denoted by $V$. By 
Theorem 6 we have $\nu_2(V) = \frac{2n}{2} = n$, as claimed.
}
\end{ex}

\begin{rema}
{\rm
($k = \mathbb C$) If $A$ is cosemisimple Hopf algebra with antipode of finite
order, it is easily seen, using Theorem 6, that 
$|\nu_2(V)| \geq 1$ for any self-dual irreducible comodule $V$.}
\end{rema}

\bigskip

\bigskip

\address{Laboratoire de Math\'ematiques Appliqu\'ees, 

Universit\'e de Pau et des Pays de l'Adour,

IPRA, Avenue de l'universit\'e, 

64000 Pau, France.}

E-mail: Julien.Bichon@univ-pau.fr

\end{document}